# Asymptotics of probability characteristics of additive arithmetic functions

Victor Volfson

ABSTRACT We study the questions of determining the asymptotics of the probabilistic characteristics of additive arithmetic functions in the paper, regardless of whether they have a limit distribution or not. Several assertions are proved about the estimation of the asymptotics of the probabilistic characteristics of strongly additive arithmetic functions, as well as additive functions of the class that have the same asymptotic behavior of the probabilistic characteristics, as for strongly additive arithmetic functions.

## 1. INTRODUCTION

Estimation of the asymptotics of arithmetic functions and moments of arithmetic functions has been and remains an urgent problem at the present time [1], [2], [3]. Moreover, the asymptotics of both additive and multiplicative arithmetic functions are investigated [4].

The method for estimating the asymptotics of the probabilistic characteristics of additive arithmetic functions having a limiting distribution was considered in [5], [6].

The questions of determining the asymptotics of the probabilistic characteristics of additive arithmetic functions, regardless of whether they have a limit distribution, are studied in this work.

An arithmetic function is a function defined on the set of natural numbers and taking values on the set of complex numbers. The name arithmetic function is due to the fact that this function expresses some arithmetic property of the natural series.

An arithmetic function is called additive if it satisfies the condition:

$$f(m) = f(p_1^{a_1}...p_t^{a_t}) = f(p_1^{a_1}) + ... + f(p_t^{a_t}) = \sum_{p^\alpha | m} f(p^\alpha)$$

A strongly additive arithmetic function is a function for which $f(p^a) = f(p)$.

If there is an arbitrary natural number $m = p_1^{a_1}...p_t^{a_t}$, then it holds for a strongly additive arithmetic function:

$$f(m) = f(p_1^{a_1}...p_t^{a_t}) = f(p_1^{a_1}) + ... + f(p_t^{a_t}) = f(p_1) + ... + f(p_t) = \sum_{p|m} f(p). \qquad (1.1)$$

Let the function $r = r(n)$ grow more slowly than any positive exponent $n$ when the value $n \to \infty$. Let's denote the mean value and variance for the function $r = r(n)$, respectively: $A(r), D(r)$.

Then it was proved [7], that for a real additive arithmetic function $f(m), m = 1,...,n$ with mean value and variance, respectively $A(n), D(n)$, under the condition $D(n) \to \infty, n \to \infty$, there is such an unboundedly increasing function $r = r(n)$ such that $\ln r(n) / \ln n \to 0$, $D(r)/D(n) \to 1$ and the limiting distribution laws:

$$P_n\{\frac{f(m)_r - A(r)}{D(r)} < x\}, P_n\{\frac{f(m) - A(n)}{D(n)} < x\}$$

is only simultaneously and in this case coincide. Following [7], we will assume that in this case $f(m)$ belongs to the class $H$ of additive arithmetic functions.

Additive arithmetic functions $f(m)$ and strongly additive $f^*(m) = \sum_{p|m} f(m)$ were investigated in [5], [6] (see 1.1), which belong to the class $H$, those have the same limiting distribution.

It is proved [7] that a strongly additive arithmetic function $f^*(m) = \sum_{p|m} f(m)$ belongs to the class $H$ if the condition is satisfied for its variance: $\ln D(n) = o(\ln \ln(n))$.

However, there are highly additive functions for which this does not work. For example, we will consider a strongly additive arithmetic function $f(m) = \sum_{p|m} \ln(p)$. In this case: $\ln D(n) \sim 2\ln \ln(n)$. Therefore, the condition: $\ln D(n) = o(\ln \ln(n))$ is not met and this strongly additive arithmetic function does not belong to the class $H$.

Although the additive arithmetic function $f(m) = \ln(m)$ and the strongly additive arithmetic function $f^*(m) = \sum_{p|m} \ln(p)$ have the same asymptotic mean and variance, but they do not have a limit distribution at all.

Let us select the class S of additive arithmetic functions for which the additive arithmetic function $f(m)$ and the corresponding strongly additive arithmetic function $f^*(m) = \sum_{p|m} f(p)$ have the same asymptotic behavior of the probabilistic characteristics.

In a particular case, the class S includes the class $H$ of additive arithmetic functions, since this condition is satisfied for it. However, the class S includes not only the class $H$ of additive arithmetic functions, but also other additive arithmetic functions, which we will consider in this work.

Formulas were obtained in [7] to determine asymptotics of the mean value of an additive arithmetic function when the value $n \to \infty$:

$$A(n) = \sum_{p^\alpha < n} \frac{f(p^\alpha)}{p^\alpha} \qquad (1.2)$$

and asymptotics of the variance of the additive arithmetic function when the value $n \to \infty$:

$$D(n) = \sum_{p^\alpha < n} \frac{f^2(p^\alpha)}{p^\alpha}. \qquad (1.3)$$

Having in mind (1.2) and (1.3), the questions of estimating the asymptotics of the probabilistic characteristics of additive arithmetic functions are reduced to finding the asymptotics of some sums of functions of primes.

2. DETERMINATION OF THE ASYMPTOTICS OF SOME SUMS OF FUNCTIONS OF PRIME NUMBERS

Let's consider the sums of prime number functions of the following form:

$$\sum_{\varphi(p) \leq x} g(p), \qquad (2.1)$$

where the function of a prime argument $\varphi(p)$ has an inverse function $\varphi^{-1}(p)$.

Then:

$$\sum_{\varphi(p)\leq x} g(p) = \sum_{p\leq \varphi^{-1}(x)} g(p). \tag{2.2}$$

Let us investigate the case (2.2), when $\varphi(p) = p^\alpha$, where $\alpha > 0$:

$$\sum_{p^\alpha \leq n, \alpha=1,\ldots,k,\ldots} g(p^\alpha) = \sum_{p\leq n} g(p) + \sum_{p\leq n^{1/2}} g(p^2) + \ldots + \sum_{p\leq n^{1/k}} g(p^k) + \ldots \tag{2.3}$$

when the value $n \to \infty$.

Having in mind (2.3), we define the asymptotic:

$$\sum_{p^\alpha \leq x, \alpha=1,\ldots,k,\ldots} \frac{f(p^\alpha)}{p^\alpha} = \sum_{p\leq n} \frac{f(p)}{p} + \sum_{p\leq n^{1/2}} \frac{f(p^2)}{p^2} + \ldots + \sum_{p\leq n^{1/k}} \frac{f(p^k)}{p^k} + \ldots \tag{2.4}$$

when the value $n \to \infty$.

We also define the asymptotic:

$$\sum_{p^\alpha \leq x, \alpha=1,\ldots,k,\ldots} \frac{f^2(p^\alpha)}{p^\alpha} = \sum_{p\leq n} \frac{f^2(p)}{p} + \sum_{p\leq n^{1/2}} \frac{f^2(p^2)}{p^2} + \ldots + \sum_{p\leq n^{1/k}} \frac{f^2(p^k)}{p^k} + \ldots, \tag{2.5}$$

when the value $n \to \infty$.

It is known [8] that the asymptotic:

$$\sum_{p\leq x} \frac{\ln p}{p} = \ln x + O(1). \tag{2.6}$$

Taking into account that the series $\sum_p \frac{\ln p^k}{p^k}$ - converges at $k \geq 2$, then based on (2.4) and (2.6) we obtain the asymptotic:

$$\sum_{p^\alpha \leq n, \alpha=1,\ldots,k,\ldots} \frac{\ln p^\alpha}{p^\alpha} = \ln n + O(1) \tag{2.7}$$

when the value $n \to \infty$.

Based on [8], the asymptotic of the expression is:

$$\sum_{p \leq n} \frac{\ln^2 p}{p} = \frac{1}{2}\ln^2 n + O(1). \qquad (2.8)$$

Then, having in mind, that the series $\sum_{p} \frac{\ln^2(p^\alpha)}{p^\alpha}$ converges at $\alpha \geq 2$, based on (2.5) and (2.8) we get:

$$\sum_{p^\alpha \leq x, \alpha=1,\ldots,k,\ldots} \frac{\ln^2(p^\alpha)}{p^\alpha} = \sum_{p \leq n} \frac{\ln^2(p)}{p} + \sum_{p \leq n^{1/2}} \frac{\ln^2(p^2)}{p^2} + \ldots + \sum_{p \leq n^{1/k}} \frac{\ln^2(p^k)}{p^k} + \ldots = \frac{1}{2}\ln^2 n + O(1) \quad (2.9)$$

when the value $n \to \infty$.

## 3. ASYMPTOTICS OF PROBALITY CHARACTERISTICS OF ADDITIVE ARITHMETIC FUNCTIONS

Assertion 1

The additive arithmetic function $f(m), m=1,\ldots,n$ and the strongly additive arithmetic function $\sum_{p|m} f(p), m=1,\ldots,n$ belong to the class $S$, i.e. have the same asymptotic behavior of the probabilistic characteristics when the following conditions are satisfied:

$$\sum_{p^\alpha \leq n, \alpha=1,\ldots,k,\ldots} \frac{f(p^\alpha)}{p^\alpha} = \sum_{p \leq n} \frac{f(p)}{p} + O(1). \qquad (3.1)$$

$$\sum_{p^\alpha \leq n, \alpha=1,\ldots,k,\ldots} \frac{f^2(p^\alpha)}{p^\alpha} = \sum_{p \leq n} \frac{f^2(p)}{p} + O(1). \qquad (3.2)$$

when the value $n \to \infty$.

Proof

Based on (1.2), we define the asymptotic of the mean value of the additive arithmetic function when the value $n \to \infty$:

$$A(n) = \sum_{p^\alpha \leq n, \alpha=1,\ldots,k,\ldots} \frac{f(p^\alpha)}{p^\alpha} = \sum_{p \leq n} \frac{f(p)}{p} + \sum_{p \leq n^{1/2}} \frac{f(p^2)}{p^2} + \ldots + \sum_{p \leq n^{1/k}} \frac{f(p^k)}{p^k} + \ldots \qquad (3.3)$$

Having in mind (1.1), we define the asymptotic of the mean value of a strongly additive arithmetic function when the value $n \to \infty$:

$$A^*(n) = \sum_{p \leq n} \frac{f(p)}{p} + O(1). \tag{3.4}$$

Taking into account (3.3) and (3.4), the equality of the asymptotics can be only when the condition (3.1) is satisfied.

Having in mind (1.3), we define the asymptotic of the variance of the additive arithmetic function when the value $n \to \infty$:

$$D(n) = \sum_{p^\alpha \leq n, \alpha=1,\ldots,k,\ldots} \frac{f^2(p^\alpha)}{p^\alpha} = \sum_{p \leq n} \frac{f^2(p)}{p} + \sum_{p \leq n^{1/2}} \frac{f^2(p^2)}{p^2} + \ldots + \sum_{p \leq n^{1/k}} \frac{f^2(p^k)}{p^k} + \ldots \tag{3.5}$$

Based on (1.1), we define the asymptotic behavior of the variance of a strongly additive arithmetic function when the value $n \to \infty$:

$$D^*(n) = \sum_{p \leq n} \frac{f^2(p)}{p} + O(1). \tag{3.6}$$

Taking into account (3.5) and (3.6), the equality of the asymptotics can be only when the condition (3.2) is satisfied.

Let us show that the additive arithmetic function $\Omega(m)$ and the strongly additive arithmetic function $\omega(m) = \sum_{p|m} \Omega(p)$ satisfy Assertion 1.

Based on (3.3), we find the asymptotic of the mean value of the additive arithmetic function when the value $n \to \infty$:

$$\sum_{p^\alpha \leq n, \alpha=1,\ldots,k,\ldots} \frac{\Omega(p^\alpha)}{p^\alpha} = \sum_{p \leq n} \frac{\Omega(p)}{p} + \sum_{p \leq n^{1/2}} \frac{\Omega(p^2)}{p^2} + \ldots + \sum_{p \leq n^{1/k}} \frac{\Omega(p^k)}{p^k} + \ldots = \sum_{p \leq n} \frac{\Omega(p)}{p} + O(1) = \sum_{p \leq n} \frac{1}{p} + O(1) = \ln \ln n + O(1). \tag{3.7}$$

Having in mind (3.4), we define the asymptotic of the mean value of a strongly additive arithmetic function $\omega(m), m = 1, \ldots, n$ when the value $n \to \infty$:

$$\sum_{p \leq n} \frac{\Omega(p)}{p} + O(1) = \sum_{p \leq n} \frac{1}{p} + O(1) = \ln \ln n + O(1). \tag{3.8}$$

Based on (3.7) and (3.8), we get:

$$\sum_{p^\alpha \le n, \alpha=1,\ldots,k,\ldots} \frac{\Omega(p^\alpha)}{p^\alpha} = \sum_{p \le n} \frac{\Omega(p)}{p} + O(1),$$

i.e. the conditions of assertion 1 are satisfied.

Using (3.5), we find the asymptotic of the variance of the additive arithmetic function when the value $n \to \infty$:

$$\sum_{p^\alpha \le n, \alpha=1,\ldots,k,\ldots} \frac{\Omega^2(p^\alpha)}{p^\alpha} = \sum_{p \le n} \frac{\Omega^2(p)}{p} + \sum_{p \le n^{1/2}} \frac{\Omega^2(p^2)}{p^2} + \ldots + \sum_{p \le n^{1/k}} \frac{\Omega^2(p^k)}{p^k} + \ldots = \sum_{p \le n} \frac{\Omega^2(p)}{p} + O(1) = \sum_{p \le n} \frac{1}{p} + O(1) = \ln \ln n + O(1). \quad (3.9)$$

Having in mind (3.6), we define the asymptotic behavior of the variance of a strongly additive arithmetic function when the value $n \to \infty$:

$$\sum_{p \le n} \frac{\Omega^2(p)}{p} + O(1) = \sum_{p \le n} \frac{1}{p} + O(1) = \ln \ln n + O(1). \quad (3.10)$$

Based on (3.9) and (3.10), we get:

$$\sum_{p^\alpha \le n, \alpha=1,\ldots,k,\ldots} \frac{\Omega^2(p^\alpha)}{p^\alpha} = \sum_{p \le n} \frac{\Omega^2(p)}{p} + O(1),$$

therefore, the conditions of assertion 1 are satisfied.

On the other hand, it is known [5] that the additive arithmetic function $\Omega(m), m = 1,\ldots,n$ and the strongly additive arithmetic function $\omega(m), m = 1,\ldots,n$ belong to the class H, therefore, they have the same limit distribution and asymptotics of the mean value and variance when the value $n \to \infty$.

Now we consider the case when an additive arithmetic function $f(m), m = 1,\ldots,n$ and a strongly additive arithmetic function $\sum_{p|m} f(p), m = 1,\ldots,n$ belong to the class $S$, but do not belong to the class H.

Assertion 2

An additive arithmetic function $f(m) = \ln(m^u)$, where $u > 0$, and a strongly additive arithmetic function $f^*(m) = \sum_{p|m} \ln(p^u)$ belong to the class $S$, but do not belong to the class H.

Proof

Taking into account that the series $\sum_{p} \frac{\ln(p^{ku})}{p^k}$ - converges at $k \geq 2$, and based on (2.7), the asymptotic of the mean value of the additive arithmetic function $f(m) = \ln(m^u), m = 1,...,n$ when the value $n \to \infty$ is equal to:

$$\sum_{p^\alpha \leq n, \alpha=1,...,k,...} \frac{\ln(p^{\alpha u})}{p^\alpha} = u \sum_{p \leq n} \frac{\ln p}{p} + O(1) = u \ln n + O(1). \qquad (3.11)$$

The asymptotic of the mean value of a strongly additive arithmetic function $f^*(m) = \sum_{p|m} \ln(p^u), m = 1,...,n$ when the value $n \to \infty$ is equal to:

$$\sum_{p \leq n} \frac{\ln(p^u)}{p} = u \sum_{p \leq n} \frac{\ln p}{p} + O(1) = u \ln n + O(1). \qquad (3.12)$$

Taking into account that the series $\sum_{p} \frac{\ln^2(p^{ku})}{p^k}$ - converges at $k \geq 2$, and based on (2.9), the asymptotic of the variance of the additive arithmetic function $f(m) = \ln(m^u), m = 1,...,n$ when the value $n \to \infty$ is equal to:

$$\sum_{p^\alpha \leq n, \alpha=1,...,k,...} \frac{\ln^2(p^{\alpha u})}{p^\alpha} = \sum_{p \leq n} \frac{\ln^2 p^u}{p} + O(1) = \frac{1}{2} u \ln^2 n + O(1). \qquad (3.13)$$

The asymptotic of the variance of a strongly additive arithmetic function $f^*(m) = \sum_{p|m} \ln(p^u), m = 1,...,n$ when the value $n \to \infty$ is equal to:

$$\sum_{p \leq n} \frac{\ln^2(p^u)}{p} = \sum_{p \leq n} \frac{\ln^2 p^u}{p} + O(1) = \frac{1}{2} u \ln^2 n + O(1). \qquad (3.14)$$

Based on (3.11), (3.12) and (3.13), (3.14), the conditions of assertion 1 are satisfied for arithmetic functions $f(m), f^*(m)$, therefore the arithmetic functions $f(m) = \ln(m^u)$, where $u > 0$, and the strongly additive arithmetic functions $f^*(m) = \sum_{p|m} \ln(p^u)$ belong to the class $S$.

Having in mind (3.14) $\ln D(n) = 2\ln\ln n + O(1) = O(\ln\ln n)$, the condition $\ln D(n) = o(\ln\ln(n))$ is not satisfied in this case and the strongly additive arithmetic function $f^*(m) = \sum_{p|m} \ln(p^u)$ does not belong to the class $H$, as well as the additive arithmetic function $f(m) = \ln(m^u)$ when a value $u > 0$ has the same value of the asymptotic of variance.

Assertion 3

Suppose that the following condition is satisfied for the additive arithmetic function $f(m), m = 1,...,n$ when the value $n \to \infty$:

$$f(m) = O(\ln(m)), \qquad (3.15)$$

then the additive arithmetic function $f(m)$ belongs to the class $S$.

Proof

Based on (3.3), the asymptotic of the mean value of the additive arithmetic function when the value $n \to \infty$ determined by the formula:

$$A(n) = \sum_{p^\alpha \leq n, \alpha=1,...,k,...} \frac{f(p^\alpha)}{p^\alpha} = \sum_{p \leq n} \frac{f(p)}{p} + \sum_{p \leq n^{1/2}} \frac{f(p^2)}{p^2} + ... + \sum_{p \leq n^{1/k}} \frac{f(p^k)}{p^k} + ...$$

The asymptotic of the mean value of the additive arithmetic function when the value $n \to \infty$ for the case (3.15) is defined as:

$$A(n) = \sum_{p^\alpha \leq n, \alpha=1,...,k,...} \frac{f(p^\alpha)}{p^\alpha} = \sum_{p \leq n} \frac{f(p)}{p} + O(1). \qquad (3.16)$$

Based on (3.5), the asymptotic of the variance of the additive arithmetic function when the value $n \to \infty$ determined by the formula:

$$D(n) = \sum_{p^\alpha \leq n, \alpha=1,...,k,...} \frac{f^2(p^\alpha)}{p^\alpha} = \sum_{p \leq n} \frac{f^2(p)}{p} + \sum_{p \leq n^{1/2}} \frac{f^2(p^2)}{p^2} + ... + \sum_{p \leq n^{1/k}} \frac{f^2(p^k)}{p^k} + ...$$

The asymptotic of the variance of the additive arithmetic function when the value $n \to \infty$ for the case (3.15) is defined as:

$$D(n) = \sum_{p^\alpha \leq n, \alpha=1,\ldots,k,\ldots} \frac{f^2(p^\alpha)}{p^\alpha} = \sum_{p \leq n} \frac{f^2(p)}{p} + O(1). \qquad (3.17)$$

Having in mind (3.16) and (3.17), the conditions of assertion 1 are satisfied for the arithmetic function $f(m), m=1,\ldots,n$ when the value $n \to \infty$, therefore, $f(m)$ belongs to the class $S$.

We consider an example of using assertion 3. Let us show that the additive arithmetic function $f(m) = \ln \varphi(m)$ satisfies the conditions of assertion 3.

It is known [9] that:

$$\varphi(m) = m \prod_{p|m}(1 - 1/p) \leq m. \qquad (3.18)$$

Having in mind (3.18):

$$\ln \varphi(m) = O(\ln(m)). \qquad (3.19)$$

Based on (3.19), the additive arithmetic function $f(m) = \ln \varphi(m)$ satisfies the conditions of assertion 3 and belongs to the class $S$, i.e. asymptotics of the probabilistic characteristics of this function coincide with asymptotics of the probabilistic characteristics of a strongly additive arithmetic function $f^*(m) = \sum_{p|m} f(p) = \sum_{p|m} \ln \varphi(p)$.

Using the formula $\varphi(m) = m \prod_{p|m}(1 - 1/p)$ we get:

$$f(p) = \ln \varphi(m) = \ln(p(1 - 1/p)) = \ln p + \ln(1 - 1/p). \qquad (3.20)$$

Having in mind (3.20) and paper [5], the asymptotic of the mean value of a strongly additive arithmetic function $f^*(m), m=1,\ldots,n$ when the value $n \to \infty$ is equal to:

$$A(n) = \sum_{p \leq n} \frac{f(p)}{p} = \sum_{p \leq n} \frac{\ln p}{p} + \sum_{p \leq n} \frac{\ln(1 - 1/p)}{p}. \qquad (3.21)$$

Taking into account [8] the asymptotic is fulfilled:

$$\sum_{p \leq n} \frac{\ln p}{p} = \ln n + O(1),$$

and the series $\sum_{p} \frac{\ln(1-1/p)}{p}$ converges, then based on (3.21) we obtain:

$$A(n) = \ln n + O(1).\tag{3.22}$$

The asymptotic of the variance of the strongly additive arithmetic function when the value $n \to \infty$ is equal to:

$$D(n) = \sum_{p \le n} \frac{f^2(p)}{p} = \sum_{p \le n} \frac{\ln^2 p}{p} + 2\sum_{p \le n} \frac{\ln p \ln(1-1/p)}{p} + \sum_{p \le n} \frac{\ln^2(1-1/p)}{p}.\tag{3.23}$$

Based on [8], the asymptotic is fulfilled $\sum_{p \le n} \frac{\ln^2 p}{p} = 0.5\ln^2 n + O(1)$, taking into account that the series $\sum_{p} \frac{\ln p \ln(1-1/p)}{p}, \sum_{p} \frac{\ln^2(1-1/p)}{p}$ - converge, having in mind (3.23) we obtain:

$$D(n) = 0.5\ln^2 n + O(1).\tag{3.24}$$

Based on assertion 3, asymptotics of the probabilistic characteristics of the additive arithmetic function $f(m) = \ln \varphi(m), m = 1,...,n$ when the value $n \to \infty$ is also determined by formulas (3.22) and (3.24) and coincides with asymptotics of the probabilistic characteristics of the additive arithmetic function $\ln m, m = 1,...,n$ when the value $n \to \infty$.

# References

1. J.Chattopadhyay, P. Darbar Mean values and moments of arithmetic functions over number fields, Research in Number Theory 5(3), 2019.

2. M. Garaev, M. Kühleitner, F. Luca, G. Nowak Asymptotic formulas for certain arithmetic functions, Mathematica Slovaca 58(3):301-308, 2008.

3. R. Brad Arithmetic functions in short intervals and the symmetric group, arXiv preprint https://arxiv.org/abs/1609.02967(2016).

4. L. Tóth, Multiplicative Arithmetic Functions of Several Variables, arXiv preprint https://arxiv.org/abs/1310.7053(2013).

5. Volfson V.L. Estimate of the asymptotic behavior of the moments of arithmetic functions having limiting normal distribution, arXiv preprint https://arxiv.org/abs/2104.10164 (2021).

6. Volfson V.L. Investigation of estimates for asymptotics of the moments of additive arithmetic functions having an arbitrary limit distribution, arXiv preprint https://arxiv.org/abs/2106.10237 (2021).

7. Kubilius I.P. Probabilistic Methods in Number Theory, Vilnius, Gospolitnauchizdat, Lithuanian SSR, 1962, 220 pp.

8. Volfson V.L. Asymptotics of summation functions, Applied Physics and Mathematics No. 3, 2020, pp. 34-37.

9. A.A. Bukhshtab. "Theory of numbers", Publishing house "Education", M., 1966, 384 p.